# Proving Probabilistic Correctness Statements: the Case of Rabin's Algorithm for Mutual Exclusion[*]

Isaac Saias[†]

Laboratory for Computer Science
Massachusetts Institute of Technology
Cambridge, MA  02139

## Abstract

The correctness of most randomized distributed algorithms is expressed by a statement of the form "some predicate of the executions holds with high probability, regardless of the order in which actions are scheduled". In this paper, we present a general methodology to prove correctness statements of such randomized algorithms. Specifically, we show how to prove such statements by a series of refinements, which terminate in a statement independent of the schedule. To demonstrate the subtlety of the issues involved in this type of analysis, we focus on Rabin's randomized distributed algorithm for mutual exclusion [6].

Surprisingly, it turns out that the algorithm does not maintain one of the requirements of the problem under a certain schedule. In particular, we give a schedule under which a set of processes can suffer lockout for arbitrary long periods.

## 1  Introduction

### 1.1  General Considerations

For many distributed system problems, it is possible to produce randomized algorithms that are better than their deterministic counterparts: they may be more efficient, have simpler structure, and even achieve correctness properties that deterministic algorithms cannot. One cost of using randomization is the increased difficulty of proving correctness of the resulting algorithms. A randomized algorithm typically involves two different types of nondeterminism – that arising from the random choices and that arising from an *adversary*. The interaction between these two kinds of nondeterminism complicates the analysis of the algorithm.

In the distributed system model considered here, each of a set of concurrent processes executes its local code and communicates with the others through a shared variable. The code can contain random choices, which leads to probabilistic branch points in the tree of executions. By assumption, the algorithm is provided at certain points of the execution with random inputs having known distributions. We can equivalently consider that *all* random choices made in a single execution are given by a parameter $\omega$ at the onset of the execution. The parameter $\omega$ thus captures the first type of nondeterminism.

For the second type, we here define the adversary $\mathcal{A}$ to be the entity controlling the order in which processes take steps. (In other work (e.g., [4]), the adversary can control other decisions, such as the contents of some messages.) An adversary $\mathcal{A}$ bases its choices on the *knowledge* it holds about the prior execution of the system. This knowledge varies according to the specifications for each given problem. In this paper, we will consider an adversary allowed to observe only certain "external" manifestations of the execution and having no access, for example, to information about local process states. We will say that an adversary is *admissible* to emphasize its specificity.

These two sources of nondeterminism, $\omega$ and $\mathcal{A}$, uniquely define an execution $\mathcal{E} = \mathcal{E}(\omega, \mathcal{A})$ of the algorithm.

Among the correctness properties one often wishes to prove for randomized algorithms are properties that state that a certain property $W$ of executions has a "high" probability of holding against all admissible adversaries. Note that the probability men-

[*]Research supported by research contracts ONR-N00014-91-J-1046, NSF-CCR-8915206 and DARPA-N00014-89-J-1988.

[†]e-mail: saias@theory.lcs.mit.edu

0

tioned in this statement is taken with respect to a probability distribution on executions. One of the major sources of complication is that there are two probability spaces that need to be considered: the space of random inputs $\omega$ and the space of random executions. Let $\mathbf{dP}$ denote the probability measure given for the space of random inputs $\omega$.

Since the evolution of the system is determined both by the (random) choices expressed by $\omega$ and also by the adversary $\mathcal{A}$, we do not have a single probability distribution on the space of all executions. Rather, for each adversary $\mathcal{A}$ there is a corresponding distribution $\mathbf{dP}_\mathcal{A}$ on the executions "compatible with" $\mathcal{A}$. High probability correctness properties of a randomized algorithm $\mathcal{C}$ are then generally stated in terms of the distributions $dP_\mathcal{A}$, in the following form. Let $W$ and $I$ be sets of executions of $\mathcal{C}$ and let $l$ be a real number in $[0, 1]$. Then $\mathcal{C}$ is correct provided that $\mathbf{P}_\mathcal{A}[W \mid I] \geq l$ for every admissible adversary $\mathcal{A}$. For a condition expressed in this form, we think of $W$ as the set of "good" (or "winning") executions, while $I$ is a set that expresses the assumptions under which the good behavior is supposed to hold.

In general, it is difficult to calculate (good bounds on) probabilities of the form $P_\mathcal{A}[W|I]$. This is because the probability that the execution is in $W$, $I$ or $W \cap I$ depends on a combination of the choices in $\omega$ and those made by the adversary. Although we assume a basic probability distribution $\mathbf{P}$ for $\omega$, the adversary's choices are determined in a more complicated way – in terms of certain kinds of knowledge of the prior execution. In particular, the adversary's choices can depend on the outcomes of prior random choices made by the processes.

The situation is much simpler in the special case where the events $W$ and $I$ are defined directly in terms of the choices in $\omega$. In this case, the desired probability can be calculated just by using the assumed probability distribution $\mathbf{dP}$.

Our general methodology for proving a high probability correctness property of the form $P_\mathcal{A}[W|I]$ consists of proving successive lower bounds:

$$\begin{aligned} \mathbf{P}_\mathcal{A}[W \mid I] &\geq \mathbf{P}_\mathcal{A}[W_1 \mid I_1] \\ &\vdots \\ &\geq \mathbf{P}_\mathcal{A}[W_r \mid I_r], \end{aligned}$$

where all the $W_i$ and $I_i$ are sets of executions, and where the last two sets, $W_r$ and $I_r$, are defined directly in terms of the choices in $\omega$. The final term, $\mathbf{P}_\mathcal{A}[W_r \mid I_r]$, is then evaluated (or bounded from below) using the distribution $\mathbf{dP}$. This methodology can be difficult to implement as it involves disentangling the ways in which the random choices made by the processes affect the choices made by the adversary.

This paper is devoted to emphasizing the need of such a rigorous methodology in correctness proofs: in the context of randomized algorithms the power of the adversary is generally hard to analyze and imprecise arguments can easily lead to incorrect statements.

As evidence supporting our point, we give an analysis of Rabin's randomized distributed algorithm [6] implementing mutual exclusion for $n$ processes using a read-modify-write primitive on a shared variable with $O(\log n)$ values. Rabin claimed that the algorithm satisfies the following correctness property: for every adversary, any process competing for entrance to the critical section succeeds with probability $\Omega(1/m)$, where $m$ is the number of competing processes. As we shall see, this property can be expressed in the general form $\mathbf{P}_\mathcal{A}[W \mid I] \geq l$. In [5], Sharir et al. gave another analysis of the algorithm, providing a formal model in terms of Markov chains; however, they did not make explicit the influence of the adversary on the probability distribution on executions.

We show that this influence is crucial: the adversary in [6] is much stronger than previously thought, and in fact, the high probability correctness result claimed in [6] does not hold.

## 1.2 Rabin's Algorithm

The problem of *mutual exclusion* [2] involves allocating an indivisible, reusable resource among $n$ competing processes. A mutual exclusion algorithm is said to guarantee *progress*[1] if it continues to allocate the resource as long as at least one process is requesting it. It guarantees *no-lockout* if every process that requests the resource eventually receives it. A mutual exclusion algorithm satisfies *bounded waiting* if there is a fixed upper bound on the number of times any competing process can be bypassed by any other process. In conjunction with the progress property, the bounded waiting property implies the no-lockout property. In 1982, Burns et al.[1] considered the mutual exclusion algorithm in a distributed setting where processes communicate through a shared read-modify-write variable. For this setting, they proved that any *deterministic* mutual exclusion algorithm that guarantees progress and bounded waiting requires that the shared variable take on at least $n$ distinct values. Shortly thereafter, Rabin published a *randomized* mutual exclusion algorithm [6] for the same shared memory distributed setting. His algorithm guarantees progress using a shared variable that takes on only $O(\log n)$ values.

It is quite easy to verify that Rabin's algorithm guarantees mutual exclusion and progress; in addition, however, Rabin claimed that his algorithm satisfies the following informally-stated *strong no-lockout property*[2].

---
[1] We give more formal definitions of these properties in Section 2.

> "If process $i$ participates in a trying round of a run of a computation by the protocol and compatible with the adversary, together with $0 \leq m-1 < n$ other processes, then the probability that $i$ enters the critical region at the end of that round is at least $c/m$, $c \sim 2/3$."                               (*)

This property says that the algorithm guarantees an approximately equal chance of success to all processes that compete at the given round. Rabin argued in [6] that a good randomized mutual exclusion algorithm should satisfy this strong no-lockout property, and in particular, that the probability of each process succeeding should depend inversely on $m$, the number of *actual* competitors at the given round. This dependence on $m$ was claimed to be an important advantage of this algorithm over another algorithm developed by Ben-Or (also described in [6]); Ben-Or's algorithm is claimed to satisfy a weaker no-lockout property in which the probability of success is approximately $c/n$, where $n$ is the total number of processes, i.e., the number of *potential* competitors.

Rabin's algorithm uses a randomly-chosen round number to conduct a competition for each round. Within each round, competing processes choose lottery numbers randomly, according to a truncated geometric distribution. One of the processes drawing the largest lottery number for the round wins. Thus, randomness is used in two ways in this algorithm: for choosing the round numbers and choosing the lottery numbers. The detailed code for this algorithm appears in Figure 1.

We begin our analysis by presenting three different formal versions of the no-lockout property. These three statements are of the form discussed in the introduction and give lower bounds on the (conditional) probability that a participating process wins the current round of competition. They differ by the nature of the events involved in the conditioning and by the values of the lower bounds.

Described in this formal style, the strong no-lockout property claimed by Rabin involves conditioning over $m$, the number of participating processes in the round. We show in Theorem 3.1 that the adversary can use this fact in a simple way to lock out any process during any round.

On the other hand, the *weak $c/n$ no-lockout property* that was claimed for Ben-Or's algorithm involves only conditioning over events that describe the knowledge of the adversary at the end of previous round. We show in Theorems 3.2 and 3.4 that the algorithm suffers from a different flaw which bars it from satisfying even this property.

We discuss here informally the meaning of this result. The idea in the design of the algorithm was to incorporate a mathematical procedure within a distributed context. This procedure allows one to select with high probability a unique random element from any set of at most $n$ elements. It does so in an efficient way using a distribution of small support ("small" means here $O(\log n)$) and is very similar to the approximate counting procedure of [3]. The mutual exclusion problem in a distributed system is also about selecting a unique element: specifically the problem is to select in *each* trying round a unique process among a set of competing processes. In order to use the mathematical procedure for this end and select a true random participating process at each round and for *all choices of the adversary*, it is necessary to discard the old values left in the local variables by previous calls of the procedure. (If not, the adversary could take advantage of the existing values.) For this, another use of randomness was designed so that, with high probability, at each new round, all the participating processes would erase their old values when taking a step.

Our results demonstrate that this use of randomness did not actually fulfill its purpose and that the adversary is able in some instances to use old lottery values and defeat the algorithm.

In Theorem 3.5 we show that the two flaws revealed by our Theorems 3.1 and 3.2 are at the center of the problem: if one restricts attention to executions where program variables are reset, and if we disallow the adversary to use the strategy revealed by Theorem 3.1 then the strong bound does hold. Our proof highlights the general difficulties encountered in our methodology when attempting to disentangle the probabilities from the influence of $\mathcal{A}$.

The algorithm of Ben-Or which is presented at the end of [6] is a modification of Rabin's algorithm that uses a shared variable of *constant* size. All the methods that we develop in the analysis of Rabin's algorithm apply to this algorithm and establish that Ben-Or's algorithm is similarly flawed and does not satisfy the $1/2en$ no-lockout property claimed for it in [6]. Actually, in this setting, the shared variables can take only two values, which allows the adversary to lock out processes with probability one, as we show

---

[2] In the statement of this property, a "trying round" refers to the interval between two successive allocations of the resource, and the "critical region" refers to the interval during which a particular process has the resource allocated to it. A "critical region" is also called a "critical section".

in Theorem 3.8.

In a recent paper [7], Kushilevitz and Rabin use our results to produce a modification of the algorithm, solving randomized mutual exclusion with $\log_2^2 n$ values. They solve the problem revealed by our Theorem 3.1 by conducting *before* round $k$ the competition that results in the control of $Crit$ by the end of round $k$. And they solve the problem revealed by our Theorem 3.2 by enforcing in the code that the program variables are reset to 0.

The remainder of this paper is organized as follows. Section 2 contains a description of the mutual exclusion problem and formal definitions of the strong and weak no-lockout properties. Section 3 contains our results about the no-lockout properties for Rabin's algorithm. It contains Theorems 3.1 and 3.2 which disprove in different ways the strong and weak no-lockout properties and Theorem 3.5 whose proof is is a model for our methodology: a careful analysis of this proof reveals exactly the origin of the flaws stated in the two previous theorems. One of the uses of randomness in the algorithm was to disallow the adversary from knowing the value of the program variables. Our Theorems 3.2 and 3.7 express that this objective is not reached and that the adversary is able to infer (partially) the value of *all* the fields of the shared variable. Theorem 3.8 deals about the simpler setting of Ben-Or's algorithm.

Some mathematical properties needed for the constructions of Section 3 are presented in an appendix (Section 4).

## 2 The Mutual Exclusion Problem

The problem of *mutual exclusion* is that of continually arbitrating the exclusive ownership of a resource among a set of competing processes. The set of competing processes is taken from a universe of size $n$ and changes with time. A solution to this problem is a distributed algorithm described by a program (code) $\mathcal{C}$ having the following properties. All involved processes run the same program $\mathcal{C}$. $\mathcal{C}$ is partitioned into four regions, $Try$, $Crit$, $Exit$, and $Rem$ which are run cyclically in this order by all processes executing $\mathcal{C}$. A process in $Crit$ is said to hold the resource. The indivisible property of the resource means that at any point of an execution, at most one process should be in $Crit$.

### 2.1 Definition of Runs, Rounds, and Adversaries

In this subsection, we define the notions of *run*, *round*, *adversary*, and *fair* adversary which we will use to define the properties of *progress* and *no-lockout*.

A *run* $\rho$ of a (partial) execution $\mathcal{E}$ is a sequence of triplets $\{(p_1, old_1, new_1), (p_2, old_2, new_2), \ldots (p_t, old_t, new_t) \ldots\}$ indicating that process $p_t$ takes the $t^{th}$ step in $\mathcal{E}$ and undergoes the region change $old_t \to new_t$ during this step (e.g., $old_t = new_t = Try$ or $old_t = Try$ and $new_t = Crit$). We say that $\mathcal{E}$ is *compatible* with $\rho$.

An admissible *adversary* for the mutual exclusion problem is a mapping $\mathcal{A}$ from the set of finite runs to the set $\{1, \ldots, n\}$ that determines which process takes its next step as a function of the current partial run. That is, the adversary is *only* allowed to see the changes of regions. For every $t$ and for every run $\rho = \{(p_1, old_1, new_1), (p_2, old_2, new_2), \ldots\}$, $\mathcal{A}[\{(p_1, old_1, new_1), \ldots, (p_t, old_t, new_t)\}] = p_{t+1}$. We then say that $\rho$ and $\mathcal{A}$ are *compatible*.

An adversary $\mathcal{A}$ is *fair* if for every execution, every process $i$ in $Try$, $Crit$, or $Exit$ is eventually provided by $\mathcal{A}$ with a step. This condition describes "normal" executions of the algorithm and says that processes can quit the competition only in $Rem$.

A *round* of an execution is the part between two successive entrances to the critical section (or before the first entrance). Formally, it is a maximal execution fragment of the given execution, containing one transition $Try \to Crit$ at the end of this fragment and no other transition $Try \to Crit$. The round of a run is defined similarly.

A process $i$ *participates* in a round if $i$ takes a step while being in its trying section $Try$.

### 2.2 The Progress and No-Lockout Properties

**Definition 2.1** An algorithm $\mathcal{C}$ that solves mutual exclusion guarantees *progress* if, for all fair adversaries, there is no infinite execution in which, from some point on, at least one process is in its $Try$ region (respectively its $Exit$ region) and no transition $Try \to Crit$ (respectively $Exit \to Rem$) occurs.

The properties that we considered thus far are non-probabilistic. The no-lockout property is probabilistic. Its formal definition requires the following notation:

Let $X$ denote any generic quantity whose value changes as the execution unfolds (e.g., a program variable). We let $X(k)$ denote the value of $X$ *just prior* to the last step ($Try \to Crit$) of the $k$th round of the execution. As a special case of this general notation, we define the following.

$\mathcal{P}(k)$ is the set of participating processes in round $k$. (Set $\mathcal{P}(k) = \emptyset$ if $\mathcal{E}$ has fewer then $k$ rounds.) The notation $\mathcal{P}(k)$ is consistent with the general notation because the set of processes participating in round $k$ is updated as round $k$ progresses: in effect the definition of this set is complete only at the end of round $k$ (this fact is at the heart of our Theorem 3.1).

$t(k)$ is the total number of steps that are taken by all the processes up to the end of round $k$.

$\mathcal{N}(k)$ is the set of executions in which all the processes $j$ participating in round $k$ reinitialize their program variables $B_j$ with a new value $\beta_j(k)$ during round $k$. ($\mathcal{N}$ stands for New-values.) $\beta_j(k)$; $k = 1, 2, \ldots$, $j = 1, \ldots, n$ is a family of iid [3] random variable whose distribution is geometric truncated at $\log_2 n + 4$ (see [6]).

For each $i$ and $k$, we let $W_i(k)$ denote the set of executions in which process $i$ enters the critical region at the end of round $k$.

We consistently use the probability theory convention according to which, for any property $\mathcal{S}$, the set of executions $\{\mathcal{E} : \mathcal{E} \text{ has property } \mathcal{S}\}$ is denoted as $\{\mathcal{S}\}$. Then:

- For each step number $t$ and each execution $\mathcal{E}$ we let $\pi_t(\mathcal{E})$ denote the run compatible with the first $t$ steps of $\mathcal{E}$. For any $t$-steps run $\rho$, $\{\pi_t = \rho\}$ represents the set of executions compatible with $\rho$. ($\{\pi_t = \rho\} = \emptyset$ if $\rho$ has fewer then $t$ steps.) We will use $\pi_k$ in place of $\pi_{t(k)}$ to simplify notation.

- Similarly, for all $m \leq n$, $\{|\mathcal{P}(k)| = m\}$ represents the set of executions having $m$ processes participating in round $k$.

The quantities $\mathcal{N}(k), \{\pi_t = \rho\}, W_i(k), \{|\mathcal{P}(k)| = m\}$ are sets of executions: for a given adversary they are random events in the probability space of random executions endowed with the measure $\mathbf{dP}_{\mathcal{A}}$.

We now present the various no-lockout properties that we want to study. A first question is to characterize relevant events $I$ over which conditioning should be done. Note first that restricting the set of executions to the ones having a certain property amounts to conditioning on this property. In particular, we will condition on the fact that process $i$ participates in round $k$. A crucial remark is that, in the worst case adversary framework that we are interested in, the adversary minimizing $\mathbf{P}_{\mathcal{A}}\left[W_i(k) \mid I\right]$ will make its choices as if "knowing" $I$. We will derive telling consequences from this fact.

We have actually in mind to compute the probability of $W_i(k)$ at different points $s_k$ of the execution.

[3] Recall that iid stands for "independent and identically distributed".

One way to go, would be to condition on the past execution. But, by our previous remark, this is tantamount to allow the adversary to this knowledge. It is then easy to see that lockout is possible. Another natural alternative that we will adopt, is to compute the probability at point $s_k$ "from the point of view of the adversary": this translates formally into conditioning over the value of the run up to point $s_k$. We will say that such a no-lockout property is *run-knowing*.

The first two definitions involve evaluating the probabilities "at the beginning of round $k$".

**Definition 2.2 (Weak, Run-knowing, Probabilistic no-lockout)** A solution to the mutual exclusion problem satisfies *weak, run-knowing* probabilistic no-lockout whenever there exists a constant $c$ such that, for every fair adversary $\mathcal{A}$, every $k \geq 1$, every $(k-1)$-round run $\rho$ compatible with $\mathcal{A}$, and every process $i$,

$$\mathbf{P}_{\mathcal{A}}\left[W_i(k) \mid \pi_{k-1} = \rho, \ i \in \mathcal{P}(k)\right] \geq c/n,$$

whenever $\mathbf{P}_{\mathcal{A}}[\pi_{k-1} = \rho, \ i \in \mathcal{P}(k)] \neq 0$.

The next property formally expresses statement (*) of Rabin. As we mentioned in our general presentation, considering rounds having $m$ participating processes corresponds to conditioning on this fact.

**Definition 2.3 (Strong, Run-knowing, Probabilistic no-lockout)** The same as in Definition 2.2 except that:

$$\mathbf{P}_{\mathcal{A}}\left[W_i(k) \mid \pi_{k-1} = \rho, \ i \in \mathcal{P}(k), \ |\mathcal{P}(k)| = m\right] \geq c/m,$$

whenever $\mathbf{P}_{\mathcal{A}}[\pi_{k-1} = \rho, \ i \in \mathcal{P}(k), \ |\mathcal{P}(k)| = m] \neq 0$.

Recalling the two interpretations of conditioning in terms of time and knowledge held by the adversary, we see that this property differs fundamentally from the previous one because, here, the adversary is provided with the number of processes due to participate in the *future* round (i.e., after $t(k-1)$). By integration over $m$, we see that an algorithm satisfying the strong property also satisfies the weak property.

The next definition is the transcription of the previous one for the case where the probability is "computed at the beginning of the execution" (i.e., $s_k = 0$ for all $k$).

**Definition 2.4 (Strong, Without knowledge, Probabilistic no-lockout)** The same as in Definition 2.2 except that:

$$\mathbf{P}_{\mathcal{A}}\left[W_i(k) \mid \ i \in \mathcal{P}(k), \ |\mathcal{P}(k)| = m\right] \geq c/m,$$

whenever $\mathbf{P}_{\mathcal{A}}[i \in \mathcal{P}(k), \ |\mathcal{P}(k)| = m] \neq 0$.

By integration over $\rho$ we see that an algorithm having the property of Definition 2.3 is stronger then one having the property of Definition 2.4. Equivalently, an adversary able to falsify Property 2.4 is stronger then one able to falsify Property 2.3.

# 3  Our Results

Here, we give a little more detail about the operation of Rabin's algorithm than we gave earlier in the introduction. At each round $k$ a new *round number* $R$ is selected at random (uniformly among 100 values). The algorithm ensures that any process $i$ that has already participated in the current round has $R_i = R$, and so passes a test that verifies this. The variable $R$ acts as an "eraser" of the past: with high probability, a newly participating process does not pass this test and consequently chooses a new random number for its *lottery value* $B_i$. The distribution used for this purpose is a geometric distribution that is truncated at $b = \log_2 n + 4$: $\mathbf{P}\bigl[\beta_j(k) = l\bigr] = 2^{-l}$ for $l \leq b-1$. The first process that checks that its lottery value is the highest obtained so far in the round, at a point when the critical section is unoccupied, takes possession of the critical section. At this point the shared variable is reinitialized and a new round begins.

The algorithm has the following two features. First, any participating process $i$ reinitializes its variable $B_i$ at most once per round. Second, the process winning the competition takes at most two steps (and at least one) after the point $f_k$ of the round at which the critical section becomes free. Equivalently, a process $i$ that takes two steps after $f_k$ and does not win the competition cannot hold the current maximal lottery value. (A process $i$ having already taken a step in round $k$ holds the current round number i.e., $R_i(k) = R(k)$. On the other hand, the semaphore $S$ is set to 0 after $f_k$. If $i$ held the highest lottery value it would pass all three tests in the code and enter the critical section.) We will take advantage of this last property in our constructions.

We are now ready to state our results. The first result states that the strong $\Omega(1/m)$ result claimed by Rabin is incorrect.

**Theorem 3.1** *The algorithm does not have the strong no-lockout property of Definition (2.4) (and hence of Definition 2.3). Indeed, there is an adversary $\mathcal{A}$ such that, for all rounds $k$, for all $m \leq n - 1$, $\mathbf{P}_{\mathcal{A}}\bigl[1 \in \mathcal{P}(k), |\mathcal{P}(k)| = m\bigr] \neq 0$ but $\mathbf{P}_{\mathcal{A}}\bigl[W_1(k) \mid 1 \in \mathcal{P}(k), |\mathcal{P}(k)| = m\bigr] = 0$.*

*Proof:* As we already remarked, the worst case adversary acts as if it knows the events on which conditioning is done. Knowing beforehand that the total number of participating processes in the round is $m$ allows the adversary to design a schedule where processes take steps in turn, where process 1 begins and where process $m$ takes possession of the critical section. Specifically, the adversary $\mathcal{A}$ does not use its knowledge about $\rho$, gives

---

**Shared variable:** $V = (S, B, R)$, where:
    $S \in \{0, 1\}$, initially 0
    $B \in \{0, 1, \ldots, \lceil \log n \rceil + 4\}$, initially 0
    $R \in \{0, 2, \ldots, 99\}$, initially *random*

**Code for $i$:**
    **Local variables:**
        $B_i \in \{0, \ldots, \lceil \log n \rceil + 4\}$, initially 1
        $R_i \in \{0, 1, \ldots, 99\}$, initially $\perp$
    **Code:**
    **while** $V \neq (0, B_i, R_i)$ **do**
        **if** $(V.R \neq R_i)$ or $(V.B < B_i)$ **then**
            $B_i \leftarrow random$
            $V.B \leftarrow max(V.B, B_i)$
            $R_i \leftarrow V.R$
        unlock; lock;
    $V \leftarrow (1, 0, random)$
    unlock;
    \* Critical Region \*\*
    lock;
    $V.S \leftarrow 0$
    $R_i \leftarrow \perp$
    $B_i \leftarrow 0$
    unlock;
    \* Remainder Region \*\*
    lock;

Figure 1: Rabin's Algorithm

---

one step to process 1 while the critical section is occupied, waits for *Exit* and then adopts the schedule $2, 2, 3, 3, \ldots, n, n, 1$. This schedule brings round $k$ to its end, because of the second property mentioned above (i.e., all processes are scheduled for two steps). For this adversary, for $2 \leq m \leq n - 1$, $|\mathcal{P}(k)| = m$ happens exactly when process $m$ wins so that $\mathbf{P}_{\mathcal{A}}[W_1(k) \cap |\mathcal{P}(k)| = m] = 0$. On the other hand, for this adversary, process $m$ wins with non zero probability, i.e., $\mathbf{P}_{\mathcal{A}}[1 \in \mathcal{P}(k) \cap |\mathcal{P}(k)| = m] \neq 0$.  ∎

The previous result is not too surprising in the light of the time interpretation given before Definition 2.2. restricting the execution to $\{|\mathcal{P}(k)| = m\}$ gives $\mathcal{A}$ too much knowledge about the *future*. We now give in Theorem 3.2 the more damaging result, stating (1) that, in spite of the randomization introduced in the round number variable $R$, the adversary is able to infer the values held in the local variables and (2) that it is able to use this knowledge to lock out a process with probability exponentially close to 1.

**Theorem 3.2** *There exists a constant $c < 1$, an adversary $\mathcal{A}$, a round $k$ and a $k - 1$-round run $\rho$ such*

that:

$$\mathbf{P}_{\mathcal{A}}[W_1(k) \mid \pi_{k-1} = \rho,\ 1 \in \mathcal{P}(k)] \leq e^{-32} + c^n.$$

We need the following definition in the proof.

**Definition 3.1** Let $l$ be a round. Assume that, during round $l$, the adversary adopts the following strategy. It first waits for the critical section to become free, then gives one step to process $j$ and then *two* steps (in any order) to $s$ other processes. (We will call these *test-processes*.) Assume that at this point the critical section is still available (so that round $l$ is not over). We then say that process $j$ is an *s-survivor* (at round $l$).

The idea behind this notion is that, by manufacturing survivors, the adversary is able to select processes having high lottery values. We now describe in more detail the selection of survivors and formalize this last fact.

In the following we will consider an adversary constructing sequentially a family of $s$-survivors for the four values $s = 2^{\log_2 n + t}$; $t = -1, \ldots, -5$. Whenever the adversary manages to select a new survivor it stores it, i.e, does not allocates it any further step until the selection of survivors is completed. ($\mathcal{A}$ actually allocates steps to selected survivors, but only very rarely, to comply with fairness. Rarely means for instance once every $nT^2$ steps, where $T$ is the expected time to select an $n/2$-survivor.) By doing so, $\mathcal{A}$ reduces the pool of test-processes still available. We assume that, at any point in the selection process, the adversary selects the test-processes *at random* among the set of processes still available. (The adversary could be more sophisticated then random, but this is not needed.) Note that a new $s$-survivor can be constructed with probability one whenever the available pool has size at least $s + 1$: it suffices to re-iterate the selection process until the selection completes successfully.

**Lemma 3.3** There is a constants $d$ such that for any $t = -5, \ldots, -1$, for any $2^{\log_2 n + t}$-survivor $j$, for any $a = 0, \ldots, 5$

$$\mathbf{P}_{\mathcal{A}}[B_j(l) = \log n + t + a] \geq d.$$

*Proof:* Let $s$ denote $\log n + t$. Let $j$ be an $s$-survivor and $i_1, i_2, \ldots, i_s$ be the test-processes used in its selection. Assume also that $j$ drew a new value $B_j(l) = \beta_j(l)$ (this happens with probability $q_1 = .99$.) Remark that $B_j(l) = \text{Max}\{B_{i_1}(l), \ldots, B_{i_s}(l), B_j(l)\}$: if this were not the case, one of the test-processes would have entered *Crit*. As the test processes are selected at random, each of them has with probability .99 a round number different from $R(l)$ and hence draws a new lottery number $\beta_j(l)$. Hence, with high probability $q_2 > 0$, 90% of them do so. The other of them keep their old lottery value $B_j(l-1)$: this value, being old, has lost in previous rounds and is therefore stochastically smaller [4] then a new value $\beta_j(l)$. (An application of Lemma 4.5 formalizes this.) Hence, with probability at least $q_1 q_2$ we have the following stochastic inequality:

$$\text{Max}\{\beta_1(l), \ldots, \beta_{s \cdot 90/100}\}$$
$$\leq_{\mathcal{L}} B_j(l) \leq_{\mathcal{L}} \text{Max}\{\beta_1(l), \ldots, \beta_{s+1}(l)\}.$$

Corollary 4.4 then shows that, for $a = 0, \ldots, 5$, with probability at least $q_1 q_2$, $\mathbf{P}_{\mathcal{A}}[B_j(l) = \log_2 s] \geq q_3$ for some constant $q_3$ ($q_3$ is close to 0.01). Hence, with probability at least $d \stackrel{\text{def}}{=} q_1 q_2 q_3$, $B_j(l)$ is equal to $\log_2 s + a$. ∎

*Proof of Theorem 3.2:* The adversary uses a preparation phase to select and store some processes having high lottery values. We will, by abuse of language, identify this phase with the round $\rho$ which corresponds to it. When this preparation phase is over, round $k$ begins.

Preparation phase $\rho$: For each of the five values $\log_2 n + t, t = -5, \ldots, -1$, $\mathcal{A}$ selects in the preparation phase many ("many" means $n/20$ for $t = -5, \ldots, -2$ and $6n/20$ for $t = -1$) $2^{\log_2 n + t}$-survivors. Let $S_1$ denote the set of all the survivors thus selected. (Note that $|S_1| = n/2$ so that we have enough processes to conduct this selection). By partitioning the set of $2^{\log_2 n - 1}$-survivors into six sets of equal size, for each of the ten values $t = -5, \ldots, 4$, $\mathcal{A}$ has then secured the existence of $n/20$ processes whose lottery value is $\log_2 n + t$ with probability bigger then $d$. (By Lemma 3.3.)

Round $k$: While the critical section is busy, $\mathcal{A}$ gives a step to each of the $n/2$ processes from the set $S_2$ that it did not select in phase $\rho$. When this is done, with probability at least $1 - 2^{-32}$ (see Corollary 4.2) the program variable $B$ holds a value bigger or equal then $\log_2 n - 5$. The adversary then waits for the critical section to become free and gives steps to the processes of $S_1$ it selected in phase $\rho$. A process in $S_2$ can win access to the critical section only if the maximum lottery value $B_{S_2} \stackrel{\text{def}}{=} \text{Max}_{j \in S_2} B_j$ of all the processes in $S_2$ is strictly less then $\log_2 n - 5$ or if no process of $S_1$ holds both the correct round number $R(k)$ and the lottery number $B_{S_2}$. This consideration gives the bound predicted in Theorem 3.2 with $c = (1 - d/100)^{1/20}$. ∎

Our proof actually demonstrates that there is an adversary that can lock out, with probability exponentially close to 1, an arbitrary set of $n/2$ processes during *some* round. With a slight improvement we can derive an adversary that will succeed in locking out (with probability exponentially close to 1) a given set $S_3$ of, for example, $n/100$ processes at *all rounds*: we just need to remark that the adversary can do without this set $S_3$ during the preparation phase $\rho$. The adversary would then alternate preparation phases $\rho_1, \rho_2, \ldots$ with rounds $k_1, k_2, \ldots$ The set $S_3$ of processes would be given steps only during rounds $k_1, k_2, \ldots$ and would be locked out at each time with probability exponentially close to 1.

---

[4] A real random variable $X$ is stochastically smaller then another one $Y$ (we write that: $X \leq_{\mathcal{L}} Y$) exactly when, for all $x \in \mathbb{R}$, $\mathbf{P}[X \geq x] \leq \mathbf{P}[Y \geq x]$. Hence, if $X \leq Y$ in the usual sense, it is also stochastically smaller.

In view of our counterexample we might think that increasing the size of the shared variable might yield a solution. For instance, if the geometric distribution used by the algorithm is truncated at the value $b = 2 \log_2 n$ instead of $\log_2 n + 4$, then the adversary is not able as before to ensure a lower bound on the probability that an $n/2$-survivor holds $b$ as its lottery value. (The probability is given by Theorem 4.1 with $x = \log n$.) Then the argument of the previous proof does not hold anymore. Nevertheless, the next theorem establishes that raising the size of the shared variable does not help as long as the size stays sub-linear. But this is exactly the theoretical result the algorithm was supposed to achieve. (Recall the $n$-lower bound of [1] in the deterministic case.) Furthermore, the remark made above applies here also: a set of processes of linear size can be locked out at each time with probability arbitrarily close to 1.

**Theorem 3.4** Suppose that we modify the algorithm so that the set of possible round numbers used has size $r$ and that the set of possible lottery numbers has size $b$ ($\log_2 n + 4 \leq b \leq n$). Then there exists positive constants $c_1$ and $c_2$, an adversary $\mathcal{A}$, and a run $\rho$ such that
$$\mathbf{P}_{\mathcal{A}}[W_1(k) \mid \pi_{k-1} = \rho, \ 1 \in \mathcal{P}(k)] \leq$$
$$e^{-32} + e^{-c_1 n/r} + c_2 \frac{r}{n^2} \ .$$

*Proof:* We consider the adversary $\mathcal{A}$ described in the proof of theorem 3.2: for $t = -5, \ldots, -2$, $\mathcal{A}$ prepares a set $T_t$ of $2^{\log_2 n + t}$-survivors, each of size $n/20$, and a set $T_{-1}$ of $2^{\log_2 n - 1}$-survivors; the size of $T_{-1}$ is $6/20 n$. (We can as before think of this set as being partitioned into six different sets.) We let $\eta$ stand for $6/20$ in the sequel.

Let $p_l$ denote the probability that process 1 holds $l$ as its lottery value after having taken a step in round $k$. For any process $j$ in $S_{-1}$ let also $q_l$ denote the probability that process $j$ holds $l$ as its lottery value at the end of the preparation phase $\rho$.

The same reasoning as in Theorem 3.2 then leads to the inequality:
$$\mathbf{P}_{\mathcal{A}}[W_1(k) \mid \pi_{k-1} = \rho, \ 1 \in \mathcal{P}(k)] \leq$$

$$e^{-32} + (1 - e^{-32})(1 - d/r)^{n/20}$$
$$+ \sum_{l \geq \log_2 n + 5} p_l (1 - \frac{q_l}{r})^{\eta n} \ .$$

Write $l = \log_2 n + x - 1 = \log_2(n/2) + x$. Then, as is seen in the proof of Corollary 4.4, $q_l = e^{-2^{1-\zeta}} 2^{1-\zeta}$ for some $\zeta \in (x, x+1)$. For $l \geq \log_2 n + 5$, $x$ is at least 6 and $e^{-2^{1-\zeta}} \sim 1$ so that $q_l \sim 2^{1-\zeta} \geq 2^{1-x}$. On the other hand $p_l = 2^{-l} = 2^{-x+1}/n$.

Define $\psi(x) \stackrel{\text{def}}{=} e^{-2^{1-x} \eta n/r}$ so that $\psi'(x) = e^{-2^{1-x} \eta n/r} 2^{1-x} \eta n/r$. Then:

$$\begin{aligned}
\sum_{l \geq \log_2 + 5} p_m (1 - \frac{q_m}{r})^{\eta n} &\leq 2/n \sum_{x \geq 6} 2^{-x} (1 - \frac{2^{1-x}}{r})^{\eta n} \\
&\leq 2/n \sum_{x \geq 6} 2^{-x} e^{-(\frac{2^{1-x}}{r} \eta n)} \\
&= 1/n \sum_{x \geq 6} 2^{1-x} e^{-(\frac{2^{1-x}}{r} \eta n)} \\
&= \frac{r}{\eta n^2} \sum_{x \geq 6} \psi'(x) \\
&\leq \frac{r}{\eta n^2} \int_5^\infty \psi'(x) dx \\
&= \frac{r}{\eta n^2} [\psi]_5^\infty \\
&= \frac{r}{\eta n^2} [1 - e^{-2^{-4} \eta n/r}] \\
&\leq \frac{r}{\eta n^2} \ .
\end{aligned}$$
∎

To simplify the notations in the sequel, we will let $i_1, \ldots, i_{|\mathcal{P}(k)|}$ denote the elements of $\mathcal{P}(k)$. And we will let $p_1, p_2, \ldots$ denote the sequence of processes taking steps in turn during round $k$: recall that a process $i$ can take several steps during the round.

The flaw of the protocol revealed in Theorem 3.2 is based on the fact that the variable $R$ does not act as an eraser of the past and that the adversary can use old values to defeat the algorithm. The flaw exhibited in Theorem 3.1 is based on the fact that, even when the old values are erased, the algorithm is sensitive to the order $p_1, p_2, \ldots$ in which participating processes are scheduled. The adversary can play on this order in two different ways. It can act on the fact that different scheduling strategies influence in different ways the size $m$ of the set $\mathcal{P}(k)$ (Strategy 1). And it can use the fact that, for a *given* number $m$ of participating processes, the mathematical distribution of the sequence $(\beta_i(k); i \in \mathcal{P}(k))$ is (a priori) sensitive to the ordering $p_1, p_2, \ldots$ (Strategy 2). The adversary of Theorem 3.1 specifically used strategy 1.

The next result shows that the two flaws exhibited in Theorems 3.1 and 3.2 are at the core of the prob-lem: the algorithm does have the strong no-lockout property when we precondition on the fact that the internal variables of the participating processes are reset to new values and when we bar the adversary from using strategy 1. We will actually prove this result for a slightly modified version of the algorithm. Recall in effect that the code given in Page 6 is optimized by making a participating process $i$ draw a new lottery number when it is detected that $V.B < B_i$. We will consider the "de-optimized" version of the code in which only the test $V.R \neq R_i$ ? causes of a new drawing to occur.

The next definition formalizes the restriction that we impose on the adversary. It says that the adversary commits itself to the value of $\mathcal{P}(k)$ at the beginning of round $k$.

**Definition 3.2** We say that an adversary is *restricted* when, for each round, it allocates a step to all participating processes (of this round) before the critical section becomes free. We will let $\mathcal{A}'$ (as opposed to $\mathcal{A}$) denote any such adversary.

We will make constant use of the notation $[n] \stackrel{\text{def}}{=} \{1, 2, \ldots, n\}$. Also, for any sequence $(a_j)_{j \in \mathbb{N}}$ we will write $a_i = \underset{j \in J}{\text{Umax}} a_j$ to mean that $i$ is the *only* index in $J$ for which $a_i = \underset{j \in J}{\text{Max}} a_j$.

**Theorem 3.5** For every process $i = 1, \ldots, n$, for every round $k \geq 1$, for every restricted adversary $\mathcal{A}'$ and for every $(k-1)$-round run $\rho$ compatible with $\mathcal{A}'$,
$\mathbf{P}_{\mathcal{A}'}[W_i(k) \mid \mathcal{N}(k), \pi_{k-1} = \rho, i \in \mathcal{P}(k), |\mathcal{P}(k)| = m]$
$\geq \frac{2}{3m}$, whenever
$\mathbf{P}_{\mathcal{A}'}[\mathcal{N}(k), \pi_{k-1} = \rho, i \in \mathcal{P}(k), |\mathcal{P}(k)| = m] \neq 0$.

*Proof:*
We first define the events $\mathcal{U}(k)$ and $\mathcal{U}'_J(k)$, where $J$ is any subset of $\{1, \ldots, n\}$:

$\mathcal{U}(k) \stackrel{\text{def}}{=} \{\exists! i \in \mathcal{P}(k) \text{ s.t. } B_i(k) = \underset{j \in \mathcal{P}(k)}{\text{Max}} B_j(k)\}$,

$\mathcal{U}'_J(k) \stackrel{\text{def}}{=} \{\exists! i \in J \text{ s.t. } \beta_i(k) = \underset{j \in J}{\text{Max}} \beta_j(k)\}$.

The main result established in [6] can formally be restated as:

$$\forall m \leq n, \ \mathbf{P}\left[\mathcal{U}'_{[m]}(k)\right] \geq 2/3. \qquad (1)$$

Following the general proof technique described in the introduction we will prove that :
$\mathbf{P}_{\mathcal{A}'}\left[\mathcal{U}(k) \mid \mathcal{N}(k), \pi_{k-1} = \rho, i \in \mathcal{P}(k), |\mathcal{P}(k)| = m\right]$
$= \mathbf{P}\left[\mathcal{U}'_m(k)\right]$, and that:
$\mathbf{P}_{\mathcal{A}'}\left[W_i(k) \mid \mathcal{N}(k), \pi_{k-1} = \rho, i \in \mathcal{P}(k), |\mathcal{P}(k)| = m, \mathcal{U}(k)\right]$
$= \mathbf{P}\left[\beta_i(k) = \underset{j \in [m]}{\text{Max}} \beta_j(k) \mid \mathcal{U}'_m(k)\right]$.

The events involved in the LHS of the two inequalities (e.g., $W_i(k)$, $\mathcal{U}(k)$, $\{|\mathcal{P}(k)| = m\}$, $\{\pi_{k-1} = \rho\}$, $\{i \in \mathcal{P}(k)\}$) depend on $\mathcal{A}'$ whereas the events involved in the RHS are pure mathematical events over which $\mathcal{A}'$ has no control.

We begin with some important remarks.

(1) By definition, the set $\mathcal{P}(k) = \{i_1, i_2, \ldots\}$ is decided by the restricted adversary $\mathcal{A}'$ at the beginning of round $k$: for a given $\mathcal{A}'$ and conditioned on $\{\pi_{k-1} = \rho\}$, the set $\mathcal{P}(k)$ is defined *deterministically*. In particular, for any $i$, $\mathbf{P}_{\mathcal{A}'}[i \in \mathcal{P}(k) \mid \pi_{k-1} = \rho]$ has value 0 or 1. Similarly, there is one value $m$ for which $\mathbf{P}_{\mathcal{A}'}[|\mathcal{P}(k)| = m \mid \pi_{k-1} = \rho] = 1$. Hence, for a given adversary $\mathcal{A}'$, if the random event $\{\mathcal{N}(k), \pi_{k-1} = \rho, i \in \mathcal{P}(k), |\mathcal{P}(k)| = m\}$ has non zero probability, it is equal to the random event $\{\mathcal{N}(k), \pi_{k-1} = \rho\} \stackrel{\text{def}}{=} I$.

(2) Recall that, in the modified version of the algorithm that we consider here, a process $i$ draws a new lottery value in round $k$ exactly when $R_i(k-1) \neq R(k)$. Hence, within $I$, the event $\mathcal{N}(k)$ is equal to $\{R_{i_1}(k-1) \neq R(k), \ldots, R_{i_m}(k-1) \neq R(k)\}$. On the other hand, by definition, the random variables (in short r.v.s) $\beta_{i_j}$; $i_j \in \mathcal{P}(k)$ are iid and independent from the r.v. $R(k)$. This proves that, (for a given $\mathcal{A}'$), conditioned on $\{\pi_{k-1} = \rho\}$, the r.v. $\mathcal{N}(k)$ is independent from all the r.v.s $\beta_{i_j}$. Note that $\mathcal{U}'_{\mathcal{P}(k)}(k)$ is defined in terms of (i.e., measurable with respect to) the $(\beta_{i_j}; i_j \in \mathcal{P}(k))$, so that $\mathcal{U}'_{\mathcal{P}(k)}(k)$ and $\mathcal{N}(k)$ are also independent.

(3) More generally, consider any r.v. $X$ defined in terms of the $(\beta_{i_j}; i_j \in \mathcal{P}(k))$: $X = f(\beta_{i_1}, \ldots, \beta_{i_m})$ for some measurable function $f$. Recall once more that the number $m$ and the *indices* $i_1, \ldots, i_m$ are determined by $\{\pi_{k-1} = \rho\}$ and $\mathcal{A}'$. The r.v.s $\beta_{i_j}$ being iid, for a fixed $\mathcal{A}'$, $X$ then depends on $\{\pi_{k-1} = \rho\}$ only through the value $m$ of $|\mathcal{P}(k)|$. Formally, this means that, conditioned on $|\mathcal{P}(k)|$, the r.v.s $X$ and $\{\pi_{k-1} = \rho\}$ are independent: $\mathbf{E}_{\mathcal{A}'}[X \mid \pi_{k-1} = \rho] = \mathbf{E}_{\mathcal{A}'}[X \mid |\mathcal{P}(k)| = m] = \mathbf{E}[f(\beta_1, \ldots, \beta_m)]$. (More precisely, this equality is valid for the value $m$ for which $\mathbf{P}_{\mathcal{A}}[\pi_{k-1} = \rho, |\mathcal{P}(k)| = m] \neq 0$.) A special consequence of this fact is that $\mathbf{P}_{\mathcal{A}'}[\mathcal{U}'_{\mathcal{P}(k)}(k) \mid \pi_{k-1} = \rho] = \mathbf{P}[\mathcal{U}'_{[m]}(k)]$.

Remark that, in $\mathcal{U}(k)$, the event $W_i(k)$ is the same as the event $\{B_i(k) = \underset{j \in \mathcal{P}(k)}{\text{Umax}} B_j(k)\}$. This justifies the first following equality. The subsequent ones are commented afterwards. Also, the set $I$ that we consider here is the one having a non zero probability described in Remark (1) above.

$\mathbf{P}_{\mathcal{A}'}[W_i(k) \mid \mathcal{U}(k), I]$
$= \mathbf{P}_{\mathcal{A}'}[B_i(k) = \underset{j \in \mathcal{P}(k)}{\text{Umax}} B_j(k) \mid \mathcal{U}(k), I]$
$= \mathbf{P}_{\mathcal{A}'}[\beta_i(k) = \underset{j \in \mathcal{P}(k)}{\text{Umax}} \beta_j(k) \mid \mathcal{U}'_{\mathcal{P}(k)}(k), I] \qquad (2)$
$= \mathbf{P}_{\mathcal{A}'}[\beta_i(k) = \underset{j \in \mathcal{P}(k)}{\text{Umax}} \beta_j(k) \mid \mathcal{U}'_{\mathcal{P}(k)}(k), \pi_{k-1} = \rho]$

Equation 2 is true because we condition on $\mathcal{N}(k)$ and because $\mathcal{U}(k) \cap \mathcal{N}(k) = \mathcal{U}'_{\mathcal{P}(k)}(k)$. Equation 3 is true because $\mathcal{N}(k)$ is independent from the r.v.s $\beta_{i_j}$ as is shown in Remark (2) above.

We then notice that the events $\{\beta_i(k) = \bigcup \max_{j \in \mathcal{P}(k)} \beta_j(k)\}$ and $\mathcal{U}'_{\mathcal{P}(k)}(k)$ (and hence their intersection) are defined in terms of the r.v.s $\beta_{i_j}$. From remark (3) above, the value of Eq. 3 depends only on $m$ and is therefore independent of $i$. Hence, for all $i$ and $j$ in $\mathcal{P}(k)$, $\mathbf{P}_{\mathcal{A}'}[W_i(k) \mid \mathcal{U}(k), I] = \mathbf{P}_{\mathcal{A}'}[W_j(k) \mid \mathcal{U}(k), I]$.

On the other hand, $\sum_{i \in \mathcal{P}(k)} \mathbf{P}_{\mathcal{A}'}[\beta_i(k) = \bigcup \max_{j \in \mathcal{P}(k)} \beta_j(k) \mid \mathcal{U}'_{\mathcal{P}(k)}(k), \pi_{k-1} = \rho] = 1$: indeed, one of the $\beta_{i_j}$ has to attain the maximum.

These last two facts imply that, $\forall i \in \mathcal{P}(k)$,

$$\mathbf{P}_{\mathcal{A}'}[W_i(k) \mid \mathcal{U}(k), I] = 1/m.$$

We now turn to the evaluation of $\mathbf{P}_{\mathcal{A}'}[\mathcal{U}(k) \mid I]$.

$$\mathbf{P}_{\mathcal{A}'}[\mathcal{U}(k) \mid I] = \mathbf{P}_{\mathcal{A}'}[\mathcal{U}'_{\mathcal{P}(k)}(k) \mid I] \quad (4)$$
$$= \mathbf{P}_{\mathcal{A}'}[\mathcal{U}'_{\mathcal{P}(k)}(k) \mid \pi_{k-1} = \rho] \quad (5)$$
$$= \mathbf{P}[\mathcal{U}'_{[m]}(k)] \geq 2/3. \quad (6)$$

Equation 4 is true because we condition on $\mathcal{N}(k)$. Eq. 5 is true because $\mathcal{U}'_{\mathcal{P}(k)}(k)$ and $\mathcal{N}(k)$ are independent (See Remark (2) above). The equality of Eq. 6 stems from Remark (3) above and the inequality from Eq. 1.

We can now finish the proof of Theorem 3.5.

$\mathbf{P}_{\mathcal{A}'}[W_i(k) \mid I]$
$\geq \mathbf{P}_{\mathcal{A}'}[W_i(k), \mathcal{U}(k) \mid I]$
$= \mathbf{P}_{\mathcal{A}'}[W_i(k) \mid \mathcal{U}(k), I] \mathbf{P}_{\mathcal{A}'}[\mathcal{U}(k) \mid I] \geq 2/3 \, m$.

∎

We discuss here the lessons brought by our results. (1) Conditioning on $\mathcal{N}(k)$ is equivalent to force the algorithm to refresh all the variables at each round. By doing this, we took care of the undesirable lingering effects of the past, exemplified in Theorems 3.2 and 3.4. (2) It is *not* true that:

$$\mathbf{P}_{\mathcal{A}}\Big[\beta_i(k) = \max_{j \in \mathcal{P}(k)} \beta_j(k) \mid \mathcal{U}'_{\mathcal{P}(k)}(k), |\mathcal{P}(k)| = m\Big] = \mathbf{P}\Big[\beta_i(k) = \max_{j \in [m]} \beta_j(k) \mid \mathcal{U}'_{[m]}(k)\Big],$$

i.e., that the adversary has no control over the event $\{\beta_i(k) = \max_{j \in \mathcal{P}(k)} \beta_j(k)\}$. (This was Rabin's statement in [6].)

Indeed, the latter probability is equal to $1/m$ whereas we proved in Theorem 3.1 that there is an adversary for which the former is 0 when $m \leq n - 1$.

The crucial remark explaining this apparent paradox is that, implicit in the expression $\mathbf{P}_{\mathcal{A}}[\beta_i(k) = \max_{j \in \mathcal{P}(k)} \beta_j(k) \mid \ldots]$, is the fact that the random variables $\beta_j(k)$ (for $j \in \mathcal{P}(k)$) are compared to each other in a specific way decided by $\mathcal{A}$, before one of them reveals itself to be the maximum. For instance, in the example constructed in the proof of Theorem 3.1, when $j$ takes a step, $\beta_j(k)$ is compared *only* to the $\beta_l(k)$; $l \leq j$, and the situation is not symmetric among the processes in $\mathcal{P}(k)$.

But, if the adversary is restricted as in our Definition 3.2, the symmetry is restored and the strong no-lockout property holds.

Rabin and Kushilevitz used these ideas from our analysis to produce their algorithm [7].

In our last Theorem 3.5 we used the restriction on the adversary $\mathcal{A}'$ mostly to derive a $1/m$ bound. If we consider a general adversary $\mathcal{A}$ it is interesting to note that we can still ensure the weak lockout-property:

**Theorem 3.6** For every process $i = 1, \ldots, n$, for every round $k \geq 1$, for every adversary $\mathcal{A}$ and for every $(k-1)$-round run $\rho$ compatible with $\mathcal{A}$,
$\mathbf{P}_{\mathcal{A}}\Big[W_i(k) \mid \mathcal{N}(k), \pi_{k-1} = \rho, i \in \mathcal{P}(k)\Big] \geq .1/n$,
whenever $\mathbf{P}_{\mathcal{A}}\Big[\mathcal{N}(k), \pi_{k-1} = \rho, i \in \mathcal{P}(k)\Big] \neq 0$.

*Proof:* Omitted. ∎

This theorem holds also if, as in the context of theorem 3.4, the algorithm uses $b$ lottery numbers. This shows that the result of Theorem 3.6 is not trivial: indeed, when $b = 2\log 2$, the probability $\mathbf{P}[\beta_i(k) = b]$ of drawing the highest possible number is a $o(1/n)$. One of the difficulties of the proof is that the apparently innocuous event $\{i \in \mathcal{P}(k)\}$ is in the future of the point $t(k-1)$ at which the probability is estimated: the adversary could conceivably also use this fact to ensure some specific values of the variables when $i$ participates.

Our Theorems 3.1, 3.2 and 3.4 explored how the adversary can gain and use knowledge of the lottery values held by the processes. The next theorem states that the adversary is similarly able to derive some knowledge about the round numbers, contradicting the claim in [6] that "because the variable $R$ is randomized just before the start of the round, we have with probability 0.99 that $R_i \neq R$." Note that, expressed in our terms, the previous claim translates into $R(k) \neq R_i(k-1)$.

**Theorem 3.7** There exists an adversary $\mathcal{A}$, a round $k$, a step number $t$, a run $\rho_t$, compatible with $\mathcal{A}$, having $t$ steps and in which round $k$ is under way such that

$$\mathbf{P}_{\mathcal{A}}[R(k) \neq R_1(k-1) \mid \pi_t = \rho_t] < .99.$$

*Proof:*

We will write $\rho_t = \rho'\rho$ where $\rho'$ is a $k-1$-round run and $\rho$ is the run fragment corresponding to the $k$th round under way. Assume that $\rho'$ indicates that, before round $k$, processes $1, 2, 3, 4$ participated *only* in round $k-1$, and that process 5 never participated before round $k$. Furthermore, assume that during round $k-1$ the following pattern happened: $\mathcal{A}$ waited for the critical region to become free, then allocated one step in turn to processes $2, 1, 1, 3, 3, 4, 4$; at this point 4 entered the critical region. (All this is indicated in $\rho'$.) Assume also that the partial run $\rho$ into round $k$ indicates that the critical region became free before any competing process was given a step, and that the adversary then allocated one step in turn to processes $5, 3, 3$, and that, after 3 took its last step, the critical section was still free. We will establish that, at this point,

$$\mathbf{P}_{\mathcal{A}}[R(k) \neq R_1(k-1) \mid \pi_t = \rho'\rho] < .99 \ .$$

By assumption $k-1$ is the last (and only) round before round $k$ where processes $1, 2, 3$ and $4$ participated. Hence $R_1(k-1) = R_2(k-1) = R_3(k-1) = R(k-1)$. To simplify the notations we will let $R'$ denote this common value. Similarly we will write $\beta'_1, \beta'_2, \ldots$ in place of $\beta_1(k-1), \beta_2(k-1), \ldots$ We will furthermore write $\beta_1, \beta_2, \ldots$ in place of $\beta_1(k), \beta_2(k), \ldots$ and $B, R$ in place of $B(k), R(k)$.

Using Bayes' rule gives us:

$$\mathbf{P}_{\mathcal{A}}[R \neq R' \mid \rho', \rho]$$
$$= \frac{\mathbf{P}_{\mathcal{A}}[R \neq R' \mid \rho'] \mathbf{P}_{\mathcal{A}}[\rho \mid \rho', R \neq R']}{\mathbf{P}_{\mathcal{A}}[\rho \mid \rho']}. \quad (7)$$

In the numerator, the first term $\mathbf{P}_{\mathcal{A}}[R \neq R' \mid \rho']$ is equal to 0.99 because $R$ is uniformly distributed and independent from $R'$ and $\rho'$. We will use this fact another time while expressing the value of $\mathbf{P}_{\mathcal{A}}[\rho \mid \rho']$:

$$\mathbf{P}_{\mathcal{A}}[\rho \mid \rho']$$
$$= \mathbf{P}_{\mathcal{A}}[\rho \mid \rho', R \neq R'] \mathbf{P}_{\mathcal{A}}[R \neq R' \mid \rho']$$
$$+ \mathbf{P}_{\mathcal{A}}[\rho \mid \rho', R = R'] \mathbf{P}_{\mathcal{A}}[R = R' \mid \rho']$$
$$= 0.99 \, \mathbf{P}_{\mathcal{A}}[\rho \mid \rho', R \neq R'] \quad (8)$$
$$+ 0.01 \, \mathbf{P}_{\mathcal{A}}[\rho \mid \rho', R = R'].$$

• Consider first the case where $R \neq R'$. Then process 3 gets a YES answer when going through the test "$(V.R \neq R_3)$ **or** $(V.B < B_3)$", and consequently chooses a new value $B_3(k) = \beta_3$. Hence

$$\mathbf{P}_{\mathcal{A}}[\rho \mid \rho', R \neq R'] = \mathbf{P}[\beta_3 < \beta_5]. \quad (9)$$

• Consider now the case $R = R'$. By hypothesis, process 5 never participated in the computation before round $k$ and hence draws a new number $B_5(k) = \beta_5$. Hence:

$$\mathbf{P}_{\mathcal{A}}[\rho \mid \rho', R = R'] =$$
$$\mathbf{P}_{\mathcal{A}}[B_3(k) < \beta_5 \mid \rho', R = R']. \quad (10)$$

As processes $1, \ldots, 4$ participated *only* in round $k-1$ up to round $k$, the knowledge provided by $\rho'$ about process 3 is exactly that, in round $k-1$, process 3 lost to process 2 along with process 1, and that process 2 lost in turn to process 4, i.e., that $\beta'_3 < \beta'_2$, $\beta'_1 < \beta'_2$ and $\beta'_2 < \beta'_4$. For the sake of notational simplicity, for the rest of this paragraph we let $X$ denote a random variable whose law is the law of $\beta'_2$ conditioned on $\{\beta'_2 > \mathrm{Max}\{\beta'_1, \beta'_3\}, \beta'_2 < \beta'_4\}$. This means for instance that, $\forall x \in \mathbb{R}$,

$$\mathbf{P}[X \geq x] = \mathbf{P}\Big[\beta'_2 \geq x \ \Big|\ \beta'_2 > \mathrm{Max}\{\beta'_1, \beta'_3\}, \ \beta'_2 < \beta'_4\Big].$$

When 3 takes its first step within round $k$, the program variable $V.B$ holds the value $\beta_5$. As a consequence, 3 chooses a new value when and exactly when $B_3(k-1)(=\beta'_3)$ is strictly bigger then $\beta_5$. (The case $\beta'_3 = \beta_5$ would lead 3 to take possession of the critical section at its first step in round $k$, in contradiction with the definition of $\rho$; and the case $\beta'_3 < \beta_5$ leads 3 to keep its "old" lottery value $B_3(k-1)$.) From this we deduce that:

$$\mathbf{P}_{\mathcal{A}}[B_3(k) < \beta_5 \mid \rho', R = R'] = \mathbf{P}[\beta'_3 < \beta_5 \mid \beta'_3 < X]$$
$$+ \mathbf{P}[\beta'_3 > \beta_5, \ \beta_3 < \beta_5 \mid \beta'_3 < X]. \quad (11)$$

Using Lemma 4.5 we derive that:

$$\mathbf{P}[\beta'_3 < \beta_5 \mid \beta'_3 < X] \geq \mathbf{P}[\beta'_3 < \beta_5].$$

On the other hand $\mathbf{P}[\beta'_3 < \beta_5] = \mathbf{P}[\beta_3 < \beta_5]$ because all the random variables $\beta_i(j), i = 1, \ldots, n, j \geq 1$ are iid. Taking into account the fact that the last term of equation 11 is non zero, we have then established that:

$$\mathbf{P}_{\mathcal{A}}[B_3(k) < \beta_5 \mid \rho', R = R'] > \mathbf{P}[\beta_3 < \beta_5]. \quad (12)$$

Combining Equations 9, 10 and 12 yields:

$$\mathbf{P}_{\mathcal{A}}[\rho \mid \rho', R = R'] > \mathbf{P}_{\mathcal{A}}[\rho \mid \rho', \ R \neq R'].$$

Equation 8 then shows that $\mathbf{P}_{\mathcal{A}}[\rho \mid \rho'] > \mathbf{P}_{\mathcal{A}}[\rho \mid \rho', \ R \neq R']$. Plugging this result into Equation 7 finishes the proof. ∎

We finish with a result showing that all the problems that we encountered in Rabin's algorithm carry over for Ben-Or's algorithm. Ben-Or's algorithm is cited at the end of [6]. The code of this algorithm is the same as the one of Rabin with the following modifications. All variables $B, R, B_i, R_i; 1 \leq i \leq n$ are boolean variables, initially 0. The distribution of the lottery numbers is also different but this is irrelevant for our discussion.

We show that Ben-Or's algorithm does not satisfy the weak no-lockout property of Definition 2.2. The situation is much simpler then in the case of Rabin's algorithm: here all the variables are boolean so that a simple reasoning can be worked out.

**Theorem 3.8 (Ben Or's Alg.)** There is an adversary $\mathcal{A}$, a step number $t$ and a run $\rho_t$ compatible with $\mathcal{A}$ such that

$$\mathbf{P}_{\mathcal{A}}\Big[W_2(k) \mid \pi_t = \rho_t,\ 2 \in \mathcal{P}(k)\Big] = 0\ .$$

*Proof:* Assume that we are in the middle of round 3, and that the run $\rho_t$ indicates that (at time 0 the critical section was free and then that) the schedule 1 2 2 3 3 was followed, that at this point 3 entered in $Crit$, that it left $Crit$, that at this point the schedule 4 1 1 5 5 was followed, that 5 entered and then left $Crit$, that 6 4 4 then took a step and that at this point $Crit$ is still free.

Without loss of generality assume that the round number $R(1)$ is 0. Then $R_2(1) = 0$, $B_1(1) = 1$ and $B_2(1) = 0$: if not 2 would have entered in $Crit$. In round 2 it then must be the case that $R(2) = 1$. Indeed if this was not the case then 1 would have entered the critical section. It must then be the case that $B_1(2) = 0$ and $B_4(2) = 1$. And then that $B_6(3) = 1$ and $R(3) = 0$: if this was not the case then 4 would have entered in $Crit$ in the 3rd round.

But at this point, 2 has *no chance* to win if scheduled to take a step! ∎

**Acknowledgments** I am deeply indebted to Nancy Lynch who suggested the problem and who constantly assisted me: this paper is hers too.

# References

[1] Burns J., Fischer M., Jackson P., Lynch N. and Peterson G. Data requirements for implementation of n- process mutual exclusion using a single shared variable. *Journal of the ACM*, 29:183-205, (1982).

[2] E. Dijkstra. Solution of a Problem in Concurrent Programming Control. *Communications of the ACM*, 321, (1966).

[3] Flajolet P. and Martin N. Probabilistic Counting Algorithms for Data Base Applications. *Journal of Computer and System Sciences*, 31:182-209, (1985).

[4] Graham R. and Yao A. On the Improbability of Reaching Byzantine Agreements *Proc. 21st ACM Symp. on Theory of Computer Science* 467-478 (1989).

[5] Hart S., Sharir M. and Pnueli A. Termination of Probabilistic Concurrent Programs *ACM Transactions on Programming Languages and Systems*, Vol 5, Num 3:356-380, (1983).

[6] Michael Rabin. N-process mutual exclusion with bounded waiting by 4 log N- shared variable. *Journal of Computation and System Sciences*, 25:66-75 (1982).

[7] Rabin M. and Kushilevitz E. Randomized Mutual Exclusion Algorithm Revisited *This proceedings*

[8] Saias I. and Lynch N. An Analysis of Rabin's Randomized Mutual Exclusion Algorithm. *MIT/LCS/TM-462* (1991).

# 4 Appendix

Theorem 4.1 and its corollaries are used in the construction of the adversary in Theorem 3.2 and Theorem 3.4. Lemma 4.5 is used mostly in the proof of Theorem 3.7. The proofs can be found in [8].

**Definition 4.1** For any sequence $(a_i)_{i \in \mathbb{N}}$ we denote $\mathsf{Max}_s a_i \stackrel{\text{def}}{=} \mathsf{Max}\{a_1, a_2, \ldots, a_s\}$.

In this section the sequence $(\beta_i)$ is a sequence of iid geometric random variables:

$$\mathbf{P}[\beta_i = l] = \frac{1}{2^l};\ l = 1, 2, \ldots$$

The following results are about the distribution of the extremal function $\mathsf{Max}_s \beta_i$. The same probabilistic results hold for iid random variables $(\beta_i')$, having the truncated distribution used by Rabin: we just need to truncate at $\log_2 n + 4$ the random variables $\beta_i$ and the values that they take. This does not affect the probabilities because, by definition, $\mathbf{P}[\beta_i'(k) = \log_2 n + 4] = \sum_{l \geq \log_2 n + 4} \mathbf{P}[\beta_i = l]$.

**Theorem 4.1** For $\frac{2^{1-x}}{s} \leq 1/2$ we have the following approximation:

$$A \stackrel{\text{def}}{=} \mathbf{P}[\mathsf{Max}_s \beta_i \geq \log_2 s + x] \sim 1 - e^{-2^{1-x}}\ .$$

$$\left| A - e^{-2^{1-x}} \right| \leq e^{-2^{1-x}} \frac{4^{1-x}}{s}\ .$$

**Corollary 4.2** $\mathbf{P}[\mathsf{Max}_s \beta_i \geq \log_2 s - 4] \geq 1 - e^{-32}$.

**Corollary 4.3** $\mathbf{P}[\mathsf{Max}_s \beta_i \geq \log_2 s + 8] \leq 0.01\ .$

**Corollary 4.4** $\mathbf{P}[\mathsf{Max}_s \beta_i = \log_2 s] \geq 0.17$,
$\mathbf{P}[\mathsf{Max}_s \beta_i = \log_2 s + l] \geq 0.01,\ \forall l = 1, \ldots, 5\ .$

**Lemma 4.5** Let $B$ and $A$ be any real-valued random variables. Then

$$\forall x \in \mathbb{R},\ \mathbf{P}[B \geq x \mid B \leq A] \leq \mathbf{P}[B \geq x].\ ^5$$

---

[5] We use the convention that $0/0 = 0$ whenever this quantity arises in the computation of conditional probabilities.